\documentclass[a4paper,11pt]{amsart}
\usepackage{amsmath}
\usepackage{amssymb}
\usepackage{amstext}
\usepackage{amsbsy}
\usepackage{amsopn}
\usepackage{amsthm}
\usepackage{amscd}
\usepackage{amsxtra}
\usepackage{amsfonts}

\usepackage{ifthen}
\usepackage{xspace}
\usepackage{fancyheadings}
\usepackage{layout}
\usepackage{hyperref}

\newcommand{\preprintserver}[2]{\href{http://xxx.lanl.gov/abs/math/#2}{#1/#2}}

\input xy
\xyoption{matrix}
\xyoption{arrow}
\xyoption{curve}
\newdir{ (}{{}*!/-5pt/@^{(}}
\newcommand{\xycenter}[1]{\begin{center}
                          \mbox{\xymatrix{#1}}
                          \end{center}
                         }

\newboolean{xlabels}
\newcommand{\xlabel}[1]{
                        \label{#1}
                        \ifthenelse{\boolean{xlabels}}
                                   {\marginpar{#1}}
                                   {}
                       }

\newcommand{\CZ}{\mathbb{C}}

\newcommand{\PZ}{\mathbb{P}}

\newcommand{\TZ}{\mathbb{T}}

\newcommand{\sO}{{\mathcal O}}

\newcommand{\from}{\leftarrow}

\newboolean{probleme}
\newcommand{\problem}[1]
           {\ifthenelse{\boolean{probleme}}
                       {{\bf(PROBLEM: #1)\bf}}
                       {}
           }
\newboolean{zukuenftiges}
\newcommand{\zukunft}[1]
           {\ifthenelse{\boolean{zukuenftiges}}
                       {{\bf(AUSBAUM\"OGLICHKEIT: #1)\bf}}
                       {}
           }
\newboolean{extras}
\newcommand{\extra}[1]
           {\ifthenelse{\boolean{extras}}
                       {{\bf EXTRA #1 EXTRA\bf}}
                       {}
           }

\DeclareMathOperator{\Img}{Im}

\DeclareMathOperator{\Tor}{Tor}
\DeclareMathOperator{\rank}{rank}

\DeclareMathOperator{\Syz}{Syz}

\DeclareMathOperator{\id}{id}

\DeclareMathOperator{\supp}{supp}

\theoremstyle{plain}
\begingroup

\newtheorem{thm}{Theorem}

\newtheorem{lem}[thm]{Lemma}

\newtheorem{prop}[thm]{Proposition}
\newtheorem{conj}[thm]{Conjecture}

\numberwithin{thm}{subsection} 
\endgroup

\newtheorem*{thm*}{Theorem}
\newtheorem*{conj*}{Conjecture}
\newtheorem*{verm*}{Vermutung}

\theoremstyle{definition}
\newtheorem{defn}[thm]{Definition}
\newtheorem{rem}[thm]{Remark}
\newtheorem{example}[thm]{Example}

\numberwithin{equation}{section}

\newcommand{\nosubsections}{\renewcommand{\thethm}{\thesection.\arabic{thm}}
                            \setcounter{thm}{0}
                           }

\newcommand{\cref}[3]{(\ref{#1}, #2 \ref{#3})}

\date{\today}

\usepackage{amssymb,amsmath}
\setboolean{probleme}{true}

\begin{document}

\title{Last syzygies of 1-generic spaces}

\address{Institiut f\"ur Mathematik\\ 
          Universit\"at Hannover\\ 
          Welfengarten 1\\ 
          D-30167 Hannnover 
         }

\email{bothmer@math.uni-hannover.de}

\urladdr{http://btm8x5.mat.uni-bayreuth.de/\textasciitilde bothmer}

\thanks{Supported by Marie Curie Fellowship HPMT-CT-2001-001238}

\author{Hans-Christian Graf v. Bothmer}

\begin{abstract}
Consider a determinantal variety $X$ of expected codimension definend by the maximal minors
of a matrix $M$ of linear forms.
Eisenbud and Popescu have conjectured that $1$-generic matrices $M$ are
characterised by the property that the syzygy ideals $I(s)$ of all last syzygies $s$ of $X$
coincide with $I_X$. In this note we prove a geometric version of this 
characterization, i.e. that $M$ is $1$-generic if and only if 
the syzygy varieties $\Syz(s)=V(I(s))$ of all last syzyzgies have the same support as $X$.
\end{abstract} 

\maketitle

\section{Introduction}
\nosubsections

Syzygies and minimal free resolutions of projective varieties where introduced to algebraic geometry by Hilbert in 1890 to define what we now call the Hilbertpolynomial. With Buchbergers algorithm \cite{Buchberger} and Schreyers algorithm \cite{SchreyerDiplom} these minimal free resolutions and their syzygyspaces 
can be calculated explicitly in many examples. General results about the dimension of the
syzygy spaces have been conjectured by Green \cite{GreenKoszul} and have resulted in the remarkable
proof of Green's Conjecture for general canonical curves by Voisin \cite{VoisinK3}, \cite{VoisinOdd}.

Much less is know about the geometric interpretation of individual syzygies. Two constructions
of geometric objects associated to a $p$-th syzygy $s$ are known, if $s$ is a syzygy in the
degree $d+1$ linear strand of a projective variety $X \subset \PZ^n$ that does not lie on
hypersurfaces of degree $d$:

\begin{enumerate}
\item the syzygy variety $\Syz(s)$ which is the vanishing locus of a twisted $(p-1)$-form
that can be associated to $s$.
\item the linear space $L(s)$ cut out by the linear forms involved in $s$. The codimension
of $L(s)$ is also called the rank of $s$
\end{enumerate}

Both $X$ and $L(s)$ are contained in $\Syz(s)$. The syzygy varieties of low rank
syzygies have been classified \cite{Sch91}, \cite{Ehbauer}, \cite{EusenSchreyer}, \cite{HCK3}. 
For example the syzygy variety of a rank $p+1$ syzygy in the degree $2$ linear
strand is
a scroll and the fibers of the scroll cut out a pencil of divisors on $X$. This connection
between pencils and divisors was one of the motivations for Green's conjecture. The syzygy
variety of a rank $p+2$-syzygy in the degree $2$ linear strand is the union of $L(s)$ with a linear section of the Grassmannian $G(2,p+3)$. The
universal subbundle on the Grassmanian restricts to a special rank $2$ bundle on $X$ with
many sections. For $K3$-Surfaces this is the Lazarsfeld-Mukai-bundle which has been used so successfully by Voisin in her proof of Green's conjecture.

In this note we consider determinantal varieties which are cut out by the
maximal minors of a matrix of linear forms $M$. Eisenbud and Popescu have studied syzygies 
of these varieties in \cite{EisenbudPopescu}. These syzygies can have arbitrary rank.
Eisenbud and Popescu consider the syzygy ideal $I(s)$ of a syzygy $s$, which can be intrinsically defined and whose vanishing locus is the syzygy variety $\Syz(s)$. They show that $M$ is $1$-generic,
if and only if all last syzygies $s$ of {\sl minimal rank} satisfy $I(s)=I_X$. They conjecture
this even holds for {\sl all} last syzygies. For $1$-generic $2\times m$-Matrices they show
that this is true. 
 
Here we show a similar geometric statement for arbitray matrices $M$, namely that $M$ is
$1$-generic if and only if  all last syzygy varieties have the same support as $X$. This is done by evaluating the
twisted $p-1$-forms associated to last syzygies at points outside of $X$ using the theory
of exterior minors of Green \cite{GreenExterior}. By calculating the tangent spaces of 
$X$ and $\Syz(s)$ in a smooth point of $X$ we can also show, that all last syzygy varieties
have the same smooth locus as $X$.

To obtain the conjecture of Eisenbud and Popescu one would have to show, that 
$\Syz(s)$ has no embedded components in the singular locus of $X$ and that the
syzygy ideal $I(s)$ is always saturated.  

The paper has two sections. In the first the definition and properties of syzygy varieties and
syzygy ideals are reviewed. The second section contains the definition and properties of $1$-generic
matrices and the proof of our theorems.

\section{Szyzgies and Syzygy Varieties}
\nosubsections

Let $X \subset \PZ^{n}$ be any projective variety over $\CZ$. 
We denote
its minimal free resolution by
\[
       \sO_X \from F_\bullet
\]
where we consider $F_\bullet$ as a bounded chain complex
\[
       F_\bullet \colon   F_0 \from F_{1} \from \dots \from 0
\]
with cohomology $\sO_X$ concentrated in degree $0$.

Since $F_\bullet$ is graded and free we can write 
\[
	F_p= \bigoplus_q F_{pq} \otimes \sO(-p-q)
\]
with $F_{pq}$ vector spaces. 
The dimensions
\[
       \beta_{pq} = \dim F_{pq}
\]
are called graded Betti numbers of $X$. Sometimes we will write more shortly
\[ 
  F_p= \bigoplus_q \sO(-p-q)^{\beta_{pq}}
\]
or collect the
graded Betti numbers $\beta_{pq}$ in a so-called Betti diagram:
\[
      \begin{array}{c|ccccc}
      &       &             &            &  &\\
      \hline
      & \beta_{00} & \beta_{10} & \beta_{20} & \dots & \\
      & \beta_{01} & \beta_{11} & \beta_{21} &       & \\
      & \vdots     &            &            &       &\\  
      &            &            &            & \beta_{pq}  & \\       
      &            &            &            & 
      \end{array}
\]
For better readability we will write a dash (``-'') if $\beta_{pq}=0$.

\begin{example}
The rational normal curve $X \subset \PZ^3$ of degree $3$ has minimal free
resolution
\[
        0 \from \sO_X \from \sO_{\PZ^3} \from \sO_{\PZ^3}(-2)^3 \to \sO_{\PZ^3}(-3)^2 \from 0.
\]
The corresponding Betti diagram is therefore
\[
    \begin{matrix}
           1 & - & - \\
           - & 3 & 2
    \end{matrix}
\]
Notice that this notation corresponds to the convention used
by the computer program Macaulay \cite{M2}.
\end{example}

\begin{rem}
By the minimality of the resolution we have
\[
       F_{pq}  = \Tor_p(S_X,\CZ)_{p+q} 
\]
where $S_X$ is the coordinate ring of $X$.
\end{rem}

\begin{rem}
The syzygy spaces $F_{pq}$ can also be calculated via Koszul cohomology
as in \cite{GreenKoszul}:
\[
     F_{pq} = K_{pq}(S_X,V) = K_{p-1,q+1}(I_X,V)
\]
where $V = H^0(\PZ^n,\sO(1))$ is the space of linear forms on $\PZ^n$
and $K_{pq}(B,V)$ is the cohomology of
the Koszulcomplex
\[
     \dots \to \Lambda^{p+1} V \otimes B_{q-1}
               \to \Lambda^{p} V \otimes B_{q}
               \to \Lambda^{p-1} V \otimes B_{q+1}
\]
for any  $S(V)$-module $B$.
\end{rem}

The last remark gives a geometric interpretation for certain syzygies:

\begin{lem}
If $X \subset \PZ^n$ is not contained in any hypersurface of degree $d$, i.e. $(I_X)_{d}=0$,
then 
\[
     F_{p,d} = H^0(\Omega_{\PZ^n}^{p-1}\otimes I_X(d+p)).
\]
In particular a $p$-th syzygy $s\in F_{p,d}$ can be interpreted as a twisted $(p-1)$-form
that vanishes on $X$.
\end{lem}

\begin{proof}
By Kozulchomology we have
\begin{align*}
     F_{p,d} &= K_{p-1,d+1}(I_X,V)\\
                   &=  \ker\Bigl( \Lambda^{p-1} V \otimes H^0(I_X(d+1))
                         \to \Lambda^{p-2} V \otimes H^0(I_X(d+2)) \Bigr)
\end{align*}
since $H^0(I_X(d))=0$. This kernel can easily be identified with 
$H^0(\Omega_{\PZ^n}^{p-1}\otimes I_X(d+p))$ by considering exterior powers of
the Euler sequence.
\end{proof}

\begin{rem} The syzygies in the above lemma are the syzygies of the first non zero row
in the Betti diagramm
\end{rem}

Often the twisted $p-1$ form associated to a syzygy does vanish on a larger variety, that
contains $X$. This leads to the following definition, introduced by Ehbauer \cite{Ehbauer}:

\begin{defn}
Let $s \in F_{pq}$ be a $p$-th syzygy of $X \subset \PZ^n$ with $(I_X)_d=0$. Then
the syzygy scheme $\Syz(s)$ of $s$ is the vanishing locus of the correponding twisted
$(p-1)$-form.
\end{defn}

The ideal of a syzygy scheme can be calculated via the exterior algebra structure
of the linear strands in the minimal free resolution of $X$:

\begin{defn}
Let $s \in F_{pd}$ be a $p$-th syzygy of $X \subset \PZ^n$ with $(I_X)_d=0$. Then
the ideal
\[
     I(s) := s \wedge \Lambda^{p-1} V^* \subset (I_X)_{d+1}
\]
is called the {\sl syzygy ideal} of $s$.
\end{defn}

\begin{prop}
$V(I(s))=\Syz(s)$
\end{prop}

\begin{proof}
Choose a basis $v_i$ of $V$ and let $v_\alpha$ be the corresponding basis
of $\Lambda^{p-1}V$. Via the inclusion
\[
        F_{pd} \hookrightarrow \Lambda^{p-1} V \otimes (I_X)_{d+1}
\]
we obtain a unique representation of $s$ as sum
\[
        s = \sum_\alpha f_\alpha \otimes v_\alpha.
\]
The $f_\alpha$ are sometimes called the polynomials involved in $s$ with respect to
the chosen basis $v_\alpha$. They generate
the syzygy ideal of $s$ by definition. On the other hand $s$ vanishes as a twisted
$(p-1)$-form if and only if all $f_\alpha$ vanish.
\end{proof}

\begin{rem}
The syzygy ideal is not necessarily reduced or even saturated.
\end{rem}

Closely related to this discussion is the notion of linear forms involved in a syzygy:

\begin{defn}
Let $s \in F_{pd}$ be a $p$-th syzygy of $X \subset \PZ^n$ with $(I_X)_d=0$. Then
\[
     s \wedge \Lambda^{p-2} V^* \otimes (I_X)_{d+1}^* \subset V
\]
is called the {\sl space of linear forms involved in $s$}. The linear space cut out by these
linear forms is called $L(s)$. The codimension of $L(s)$ in $\PZ(V)$ is called the {\sl rank} of
$s$.
\end{defn}

\begin{rem}
We only mention this definition for reference. Since the methods of this paper work
for syzyzgies of arbitrary rank, we will not need it here.
\end{rem}

\section{Triple Tensors and $1$-generic Matrices of Linear Forms} \label{tripletensors}
\nosubsections

Let $A$, $B$ and $C$ be finite dimensional vector spaces of
dimensions $a$, $b$ and $c$ together with a linear map
\[
          \gamma \colon  A \otimes B \to C.
\]
$\gamma$ can be interpreted as a triple tensor 
$\gamma \in A^*\otimes B^*\otimes C$ or after choosing bases
as an $a \times b$-matrix of linear forms on $\PZ(C)$. Here we adhere to the 
Grothendiek convention of interpreting elements of $\PZ(C)$ as linear forms on $C$ or
equivalently the elements of $C$ as linear forms on $\PZ(C)$.

\begin{defn}
A non zero linear map $\CZ \to A$ is called a {\sl generalized row index} of $\gamma$
since it induces a map
\[
       \CZ \otimes B \to C
\]
which can be interpreted, up to a constant factor, as a $1 \times b$ row vector of linear-forms.

If $\CZ \to A$ is such a generalized row index, the image of $\CZ$ in $A$ under this map is a line.
We will call these images {\sl generalized rows}. 
The generalized rows form a
projective space $\PZ(A^*)$ which we call 
the {\sl row space} of $\gamma$. Similarly $\PZ(B^*)$ is the 
{\sl column space} of $\gamma$. 
\end{defn}

On the row space $\PZ(A^*)$ the triple tensor $\gamma$ induces 
a map of vector bundles
\[
           \gamma_A \colon   \sO_{\PZ(A^*)}(-1) \otimes B \to C
\]
by composing $\gamma$ with the first map of the twisted Euler sequence
\[
          0 \to \sO_{\PZ(A^*)}(-1)\otimes B 
            \to A\otimes B 
            \to \TZ_{\PZ(A^*)}(-1)\otimes B 
            \to 0
\]
on $\PZ(A^*)$.
Similarly we have
\[
           \gamma_B \colon    A \otimes  \sO_{\PZ(B^*)}(-1) \to C
\]
on the column space $\PZ(B^*)$.
From now on we will restrict our discussion to the row space $\PZ(A^*)$,
leaving the analogous constructions for the column space $\PZ(B^*)$ to
the reader.

Given a generalized row $\alpha \in \PZ(A^*)$ the restriction of
$\gamma_A$ to $\alpha$
\[
           \gamma_\alpha \colon B \to C
\]
is a map of vector spaces.

\begin{defn}
The {\sl rank} of a generalized row $\alpha$ is defined as 
$\rank \alpha := \rank \gamma_\alpha$. The image
$\Img(\gamma_\alpha) \subset C$ is called the 
{\sl space of linear forms on $\PZ(C)$ involved in $\alpha$}.
\end{defn}

\newcommand{\ymin}{Y_{min}}

\begin{rem}
The determinantal varieties
associated to $\gamma_A$ stratify the row space $\PZ(A^*)$ according
to the rank of the rows. In particular the minimal rank rows form
a closed subscheme $\ymin \subset \PZ(A^*)$. 
\end{rem}

\begin{rem}
In practice $\ymin$ is often not of expected codimension.
\end{rem} 

\begin{example} \xlabel{e-twobythree}
Consider vector spaces $A$, $B$ and $C$ of dimension $2$, $3$ and 
$4$ with basis $a_i$, $b_i$ and $c_i$. The triple tensor
\[
       \gamma \colon A \otimes B \to C
\]
with 
\begin{align*}
      \gamma(a_1 \otimes b_1) &=c_1 
      & \gamma(a_1 \otimes b_2) &=c_2 
      & \gamma(a_1 \otimes b_3) &=c_3 \\
        \gamma(a_2 \otimes b_1) &=c_2 
      & \gamma(a_2 \otimes b_2) &=c_3 
      & \gamma(a_2 \otimes b_3) &=c_4 \\
\end{align*}
can be 
represented by the matrix
\[
          \begin{pmatrix}
                c_1 & c_2 & c_3 \\
                c_2 & c_3 & c_4 \\
          \end{pmatrix}.
\]
In this basis we see two rows of rank $3$. 
Generalized rows are linear combinations
of those two. The map 
\[
       \gamma_A \colon \sO_{\PZ(A^*)}(-1) \otimes B \to C
\]
can be represented by the matrix
\[
      \begin{pmatrix}
             a_1 & a_2 & 0   & 0   \\
             0   & a_1 & a_2 & 0   \\
             0   & 0   & a_1 & a_2
      \end{pmatrix}.
\]
Since this matrix has full rank everywhere on the row-space $\PZ(A^*)$ we
see that all generalized rows of $\gamma$ have the same rank $3$.
\end{example}

\begin{defn}[$1$-generic spaces]
A triple tensor
\[
        \gamma \colon A \otimes B \to C
\]
is called {\sl $1$-generic}, if all generalized rows have rank $b$
and all generalized columns have rank $a$.
\end{defn}

\begin{example}
The $2 \times 3$-matrix considered above is $1$-generic.
\end{example}

$1$-generic triple tensors have many interesting properties. See Eisenbud
\cite{linearsectionsofdeterminantal} for a discussion. 
For example if $a \ge b$ the third induced map
\[
    \gamma_C \colon A \otimes \sO_{\PZ(C)}(-1) \to B^*
\]
on $\PZ(C)$ drops rank in expected codimension and the resulting determinantal
locus $X=X_{b-1}(\gamma_C)$ is resolved by the
Eagon-Northcott complex
\[
       I_X \from E_\bullet
\]
with $E_\bullet$ linear of degree $b$  and
\[
   E_i = E_{ib} \otimes \sO(-i-b), \quad \quad
   E_{ib} = \Lambda^{b+i} A \otimes \Lambda^{b} B \otimes S_i B.
\]
Notice that $F_{p,b-1}=E_{p-1,b}$, since here we resolve the ideal of $X_{b-1}$.

Green has observed, that the exterior minors of a $1$-generic matrix
also behave nicely:

\begin{defn}
Consider the natural map
\xycenter{
          \Lambda^n A \otimes S_n B \ar@{ (->}[r] \ar@/_1pc/[rr]_-{e_n}
          &\Lambda^n (A\otimes B) \ar[r]
          &\Lambda^n C
                 }
obtained by taking the $n$th exterior power of $\gamma$. Then the elements in 
the image of $e_n$ are called degree $n$ exterior minors of $\gamma$.
\end{defn}


\begin{prop}[Green] \label{exterior}
If $\gamma$ is $1$-generic, then $e_a$ is injective.
\end{prop}

\begin{proof}
\cite[Proposition 1.2]{GreenExterior}
\end{proof}

Lets now consider the last syzygies of $X_{b-1}$:

\begin{defn}
The syzygies $s \in E_{a-b,b}$ are called last syzygies of $X_{b-1}$ since
$E_{a-b,b}$ is the last nonzero syzygy space of $X_{b-1}$.
\end{defn}

The representation of a last syzygy of $X_{b-1}$ as an element in the corresponding
Koszul cohomoloy group can be given explicitly:

\begin{lem}[Eisenbud,Popescu]
The inclusion
\[ 
      E_{a-b,b} \hookrightarrow \Lambda^{a-b} C \otimes (I_X)_b
\]
is given by the composition

\xycenter{
       E_{a-b,b} \ar@{=}[r] 
       & \Lambda^a A \otimes \Lambda^b B \otimes S_{a-b} B 
           \ar@{ (->}[d]
       \\
        & \Lambda^b A \otimes \Lambda^b B \otimes  \Lambda^{a-b} A \otimes S_{a-b} B
            \ar[d]^{\id\otimes e_{a-b}}
        \\
         & \Lambda^b A \otimes \Lambda^b B \otimes \Lambda^{a-b} C
             \ar@{=}[r]                
         & (I_X)_b \otimes \Lambda^{a-b} C.
         }
         
\end{lem}

\begin{proof}
\cite[Theorem 2.1 and proof of Theorem 3.1]{EisenbudPopescu}
\end{proof}

Using this Eisenbud and Popescu prove

\begin{thm}[Eisenbud, Popescu]
Let $X_{b-1}(\gamma_C) \subset \PZ(C)$ be a determinantal variety of expected dimension defined by a linear map
\[
     \gamma \colon A \otimes B \to C.
\]
If for all last syzygies $s \in E_{a-b,b}$ we have $I(s) = I_X$ then $\gamma$ is
$1$-generic. 
\end{thm}

They conjecture the converse

\begin{conj}[Eisenbud, Popescu]
Let $\gamma \colon A \otimes B \to C$ be a $1$-generic map of vector spaces.
Then all last syzygies $s$ of the corresponding determinantal variety
$X := X_{b-1} \subset \PZ^n$ have the same syzygy ideal $I(s)=I_X$.
\end{conj}

which they can prove for $b=2$. Evaluating last syzygies directly, we obtain
a geometric version of this conjecture:

\begin{thm} \label{samesupport}
Let $\gamma \colon A \otimes B \to C$ be a $1$-generic map of vectorspaces.
Then all syzygy varieties $\Syz(s)$ of last syzygies $s$ of the corresponding determinantal variety
$X := X_{b-1} \subset \PZ^n$ have the same support $\supp \Syz(s) = \supp X$.
\end{thm}

\begin{proof}
Let $x \in \PZ^n$ a point not contained in $X_{b-1}$ and $s$ any last syzygy. We have
to prove, that $s$ does not vanish in $x$. 

Since $x$ is not in $X_{b-1}$ the map $\gamma_C$ has full rank in $x$. Therfore we
can choose bases of $A$, $B$ and $C$ such that $\gamma_C$ can be represented
by a matrix of linear forms 
\[
   M = \begin{pmatrix} 
            c_{11} & \dots & c_{1b} \\
            \vdots &           & \vdots \\ 
            c_{a1} & \dots & c_{ab} \\
            \end{pmatrix}
\]
such that 
\[
      M(x) = \begin{pmatrix}
                       1 & \dots & 0 \\
                        \vdots & \ddots & \vdots \\
                        0 & \dots & 1 \\
                        \hline
                        0 & \dots & 0 \\
                        \vdots & & \vdots \\
                        0 & \dots & 0   
                  \end{pmatrix}.
\]
Now by the description of Eisenbund and Popescu above, 
the representation of a last syzygy $s$ in this basis is
\[
         s = \sum_{|\beta|=b} f_\beta \otimes g_{\bar{\beta}}
\]
where $f_\beta$ is the $b \times b$-minor involving the rows $\beta_1, \dots, \beta_b$
of $M$ and $g_{\bar{\beta}}$ is a degree $a-b$ exterior minor of the remaining $(a-b)\times b$ matrix.
At $x$ all minors of $M$ except $f_{1,2,\dots,b}(x))=1$ vanish, and
\[
    s(x) = g_{b+1,\dots,a}.
\]
Since $g_{b+1,\dots,a}$ is a degree $a-b$ exterior minor of a $1$-generic $(a-b) \times b$
Matrix it is nonzero by Green's proposition \ref{exterior} above.
\end{proof}

Interestingly the converse of this theorem is also true, strenghthening the
theorem of Eisenbud and Popescu.

\begin{thm}
Let $\gamma \colon A \otimes B \to C$ be a map of vector spaces,
such that the corresponding determinantal variety $X := X_{b-1}$ has
expected codimension and for every last syzygy $s$ of $X$
we have $\supp \Syz(s) = \supp X$. Then $\gamma$ is $1$-generic.
\end{thm}

\begin{proof}
Suppose $\gamma$ is not $1$-generic. Then there exists a generalized row $\alpha$
of rank at most $b-1$. Let $A' \subset A$ be a $a-1$ dimensional
subspace that does not contain $\alpha$, and consider the induced map
\[
      \gamma' \colon A' \otimes B \to C.
\]
If the associated $\gamma'_C$ was of submaximal rank everywhere, $\gamma'$
would have a generalized column of zeros. The corresponding generalized column
of $\gamma_C$ would then be of maximal rank $1$, i.e. the vanishing set of this
column is at most of codimension $1$. Since $\supp X_{b-1}$ contains the vanishing
set of each generalised column, this would imply that $X_{b-1}$ has a component of 
codimension at most $1$,
which contradicts our assumption that $X_{b-1}$ has expected codimension.
Consequently there is a point $x \in \PZ^n$ outside of $X_{b-1}$ such that each matrix $M'$ representing
$\gamma'_C$ has full rank in $x$. 

We can therfore choose bases of $A$, $B$ and $C$ such that 
\begin{enumerate}
\item  \label{notonegeneric}  $M$ has the form
\[
        M = \begin{pmatrix}
                    0 & * & \dots & * \\ 
                    *  & * & \dots & * \\ 
                     \vdots & \vdots &  & \vdots \\ 
                    * & * & \dots & * \\ 
               \end{pmatrix}
\]
\item $M'$ consists of the last $a-1$ rows of $M$
\item$M(x)$ has the form
\[
      M(x) = \begin{pmatrix}
                        0 & \dots & 0 \\
                        \vdots & & \vdots \\
                        0 & \dots & 0 \\  
                       1 & \dots & 0 \\
                        \vdots & \ddots & \vdots \\
                        0 & \dots & 1 \\
                   \end{pmatrix}.
\]
\end{enumerate}
Now consider the syzygy $s=(b_1)^{a-b}$ where $b_1$ is the basiselement of $B$
corresponding to the first column. When we evaluate $s$ at $x$ we obtain
$s(x) = f_{a-b+1\dots a}(x) \otimes g_{1\dots a-b,s} = g_{1\dots a-b,s}$ since 
$f_{a-b+1\dots a}(x) = 1$ is the only nonzero maximal minor of $M(x)$.

The exterior minor of the upper $(a-b)\times b$
submatrix corresponding to $s=(b_1)^{a-b}$ is the wedge product of the
first $a-b$ linearforms in the first column of $M$. This wedge product vanishes
since the first of these linear forms is identically zero by property (\ref{notonegeneric}). So
$s$ is a syzygy whose syzygy variety has support outside of $X_{b-1}$.
\end{proof} 

Our methods also allow us to describe the smooth locus of all syzygy varieties:

\begin{thm}
Let $\gamma \colon A \otimes B \to C$ be a $1$-generic map of vector spaces.
Then all syzygy varieties $\Syz(s)$ of last syzygies $s$ of the corresponding determinantal variety
$X := X_{b-1} \subset \PZ^n$ have the same smooth locus, which is also the
smooth locus of  $X$.
\end{thm}

\begin{proof}
Let $s$ be any last syzygy of $X$. Since $I(s) \subset I_X$ by definition and $\supp X$ = $\supp \Syz(s)$ by theorem \ref{samesupport}, we know that the smooth locus of $\Syz(s) = V(I(s))$ is contained in the smooth locus of $X$.

For the converse Let $x \in \PZ^n$ be a point contained in the smooth locus of $X$. 
We have to prove, that the tangent space of $\Syz(s)$ in $x$ is the same as 
the tangent space of $X$ in $x$.

Since $x$ is in the smooth locus of $X_{b-1}$ the morphism  $\gamma_C$ has rank $b-1$ in $x$. Therefore we
can choose bases of $A$, $B$ and $C$ such that $\gamma_C$ can be represented
by a matrix of linear forms 
\[
   M = \begin{pmatrix} 
            c_{11} & \dots & c_{1b} \\
            \vdots &           & \vdots \\ 
            c_{a1} & \dots & c_{ab} \\
            \end{pmatrix}
\]
such that 
\[
      M(x) = \begin{pmatrix}
                       1 & \dots & 0  & 0\\
                        \vdots & \ddots & \vdots & \vdots\\
                        0 & \dots & 1 & 0\\
                        \hline
                        0 & \dots & 0 & 0\\
                        \vdots & & \vdots & \vdots\\
                        0 & \dots & 0  & 0
                  \end{pmatrix}.
\]
Now suppose $x+\epsilon y$ is a tangent vector of $X$ at $x$. Then all maximal minors
of $M(x+\epsilon y)$ have to vanish, in particular those that contain the first $b-1$ rows
and the $i$-th row ($i\ge b$):
\[
    0 = \det \begin{pmatrix}
                       1 + \epsilon c_{11}(y)& \dots & \epsilon c_{1,b-1}(y) & \epsilon c_{1b}(y)\\
                        \vdots & \ddots & \vdots & \vdots\\
                        \epsilon c_{b-1,1}(y) & \dots & 1 + \epsilon c_{b-1,b-1}(y) & \epsilon c_{b-1,b}(y)\\
                        \epsilon c_{i,1}(y) & \dots & \epsilon c_{i,1b-1}(y) & \epsilon c_{ib}(y)\\
                   \end{pmatrix}
         = \epsilon c_{ib}(y).
\]
All other minors vanish since every term of the corresponding determinant involves at least $\epsilon^2$. So $x+\epsilon y$ is tangent to $X$ if and only if $c_{bb}(y)=\dots=c_{ab}(y)=0$.

Now assume that $x+\epsilon y$ is not a tangent vector of $X$. Then we can assume after a
another base chance of $C$, that $M(x+\epsilon y)$ has the form
\[
      M(x+\epsilon y) =    \begin{pmatrix}
                       1 + \epsilon c_{11}(y)& \dots & \epsilon c_{1,b-1}(y) & 0\\
                        \vdots & \ddots & \vdots & \vdots\\
                        \epsilon c_{b-1,1}(y) & \dots & 1 + \epsilon c_{b-1,b-1}(y) & 0\\
                        \epsilon c_{b,1}(y) & \dots & \epsilon c_{b,b-1}(y) & \epsilon \\
                        \epsilon c_{b+1,1}(y) & \dots & \epsilon c_{b+1,b-1}(y) & 0 \\
                        \vdots & \ddots & \vdots & \vdots\\
                        \epsilon c_{a1}(y) & \dots & \epsilon c_{a,b-1}(y) & 0 \\
                   \end{pmatrix}
\]
As before the representation of a last syzygy $s$ in this basis is
\[
         s = \sum_{|\beta|=b} f_\beta \otimes g_{\bar{\beta}}
\]
where $f_\beta$ is the $b \times b$-minor involving the rows $\beta_1, \dots, \beta_b$
of $M$ and $g_{\bar{\beta}}$ is a degree $a-b$ exterior minor of the remaining $a-b\times b$ matrix.
At $x+\epsilon y$ all minors of $M$ except $f_{1,2,\dots,b}(x)=\epsilon$ vanish, and
\[
    s(x) = \epsilon g_{b+1,\dots,a}.
\]
Since $g_{b+1,\dots,a}$ is again a degree $a-b$ exterior minor of a $1$-generic $(a-b) \times b$
matrix it is nonzero by Green's  proposition \ref{exterior} above. Therfore $x+\epsilon y$ is not a
tangent vector of $\Syz(s)$. This shows that the tangent space of $\Syz(s)$ at $x$ is contained
in the tangent space of $X$ at $x$. Since on the other hand $\Syz(s)$ contains $X$
as scheme both tangent spaces have to coincide.
\end{proof}


\end{document}